\documentclass[a4paper,14pt]{article}

\usepackage{multicol}
\usepackage{amsthm}

\usepackage[dvips]{graphicx,color}
\usepackage{amssymb}
\usepackage{amsmath}

\begin{document}

\fontsize{14pt}{16.5pt}\selectfont

\begin{center}
\bf{Existence of a space filled with an arbitrary finite number of mutually disjoint self-similar spaces}
\end{center}

\fontsize{12pt}{11pt}\selectfont
\begin{center}
Akihiko Kitada$^{1}$, Shousuke Ohmori$^{2}$, Tomoyuki Yamamoto$^{1,2,*}$\\ 
\end{center}

\noindent
$^1$\it{Institute of Condensed-Matter Science, Comprehensive Resaerch Organization, Waseda University,
3-4-1 Okubo, Shinjuku-ku, Tokyo 169-8555, Japan}\\
$^2$\it{Faculty of Science and Engineering, Waseda~University, 3-4-1 Okubo, Shinjuku-ku, Tokyo 169-8555, Japan}\\
*corresponding author: tymmt@waseda.jp\\
~~\\
\rm
\fontsize{10pt}{11pt}\selectfont\noindent

\noindent
{\bf Abstract}\\
We discuss a sufficient condition for a space to be filled with an arbitrary finite number of self-similar spaces using a topological concept.\\

It is a great issue in the fields of materials science and geology, which are closely related to crystallography, that a space can be filled with an arbitrary finite number of grains, each of which is characterized as self-similar.
Mathematical methods for the crystallography have long been developed mainly by the group theory, because it is quite useful to discuss symmetry and periodicity in the crystal.
However, it is not necessary for the above problem to consider a periodicity of the internal elements in each grain using a concept of space group. 
We have studied the mathematical structures of self-similar condensed-matters using a point set topology \cite{kitada-csf}.
To our knowledge, such fundamental topological approach has not been applied for the study on the structures of aggregates of grains.
Here we can also discuss the tiling problems according to a same procedure for the current problem on aggregate of grains by replacing an aggregate of grains with a tiling and a grain with a tile, respectively.
The unit of the discussion is now single crystal or grain, whose characters can be defined topologically.
Therefore, it is also possible to discuss the structure of aggregates of the noncrystalline or amorphous grains using a topological method. 

Although structures of self-similar materials are often discussed using a fractal dimension in fractal science, such discussions are done only within a real space as it can be seen. 
In the present report, we employ somewhat indirect method of observation of the structures of material. 
That is, the structures are expressed indirectly through the mathematical observations of the formation of a set of equivalence classes. (see, for example, a literature by Fern$\acute{a}$ndez \cite{Fernandez}, concerning the use of quotient space in statistical physics) 
The concept of the equivalence class is familiar in the diffraction crystallography as a reciprocal lattice point \cite{Cassels}. 
Group of the lattice plane is gathered as a concept of equivalence class, and then the geometrical structure in a real space can be determined by diffraction patterns in a reciprocal space. 
The structures of materials are observed here mainly from the viewpoint of this equivalence class, that is, a decomposition space, which is similar idea as that for the diffraction analysis.
From a viewpoint of topology, any pattern can be obtained from the space which is characterized, in principle, as 0-dim, perfect, compact T$_2$ space. In the followings, we discuss this sufficient condition for a space to be filled with the self-similar spaces from the viewpoint of topological aspects.

\begin{figure}[h!]
\begin{center}
\vspace{5mm}
\includegraphics[scale=0.5]{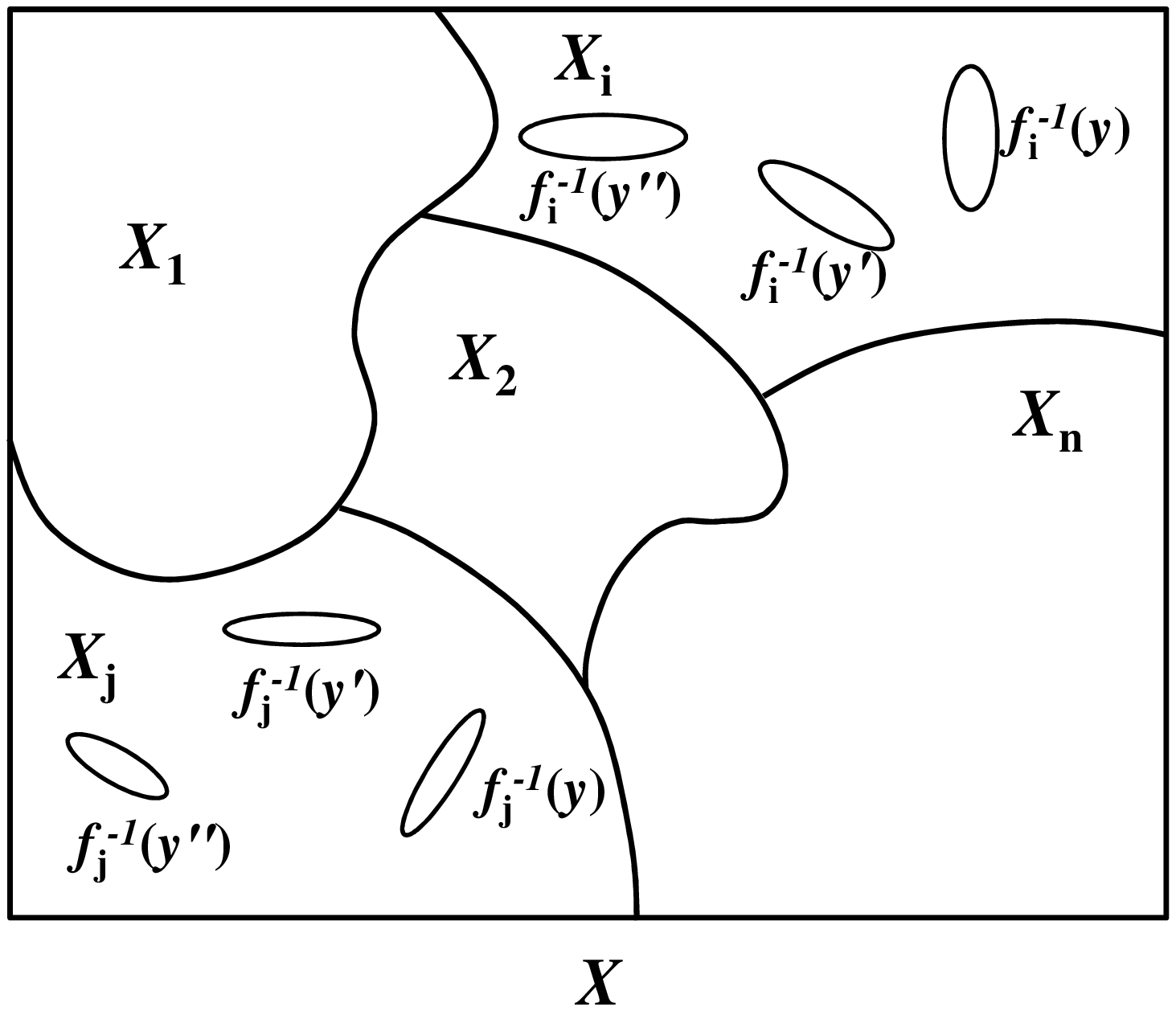}
\end{center}
\caption{A schematic explanation of the mathematical procedure. $X$ and each $X_i$ correspond to an aggregate of grains and  each grain, respectively. Equivalence class $f^{-1}_i(y)$ is a point of a self-similar decomposition space $\mathcal{D}_{f_i}$ of $X_i$. $f^{-1}_i(y)$ is not a point of $\mathcal{D}_{f_j}, i\not = j$.}
\end{figure}

Let us start the mathematical observations of the structure of material for the case shown in Fig. 1.
This figure shows the structure of aggregates of grains ($X$), where each $X_i$ corresponds to each grain.
Let $(X, \tau)=(\{0,1\}^\Lambda ,\tau_0^\Lambda ),~Card \Lambda \succ \aleph _0$
(aleph zero) be the $\Lambda -$product space of $(\{0,1\},\tau_0)$ where $\tau_0$ is a 
discrete topology for $\{0,1\}$. The topological space $(X,\tau)$ need not to be metrizable, 
that is, the relation $Card \Lambda \succ \aleph $ (aleph) may hold. In the present report, assuming the element of the grain to be a map $x:\Lambda \rightarrow \{0,1\}$ (for example, $\Lambda =\bf{N}, \bf{R},\dots$), we 
will mathematically confirm the existence of a partition $\{X_1,\dots,X_n\}$ of $X$, $X_i \in 
{(\tau\cap \Im)}-\{\phi\}$ \cite{foot2}, a decomposition space (i.e, a space of equivalence 
classes) of each $X_i$ of which is self-similar. The partition of $X$ can be regarded as a kind of aggregate of grains in an abstract sense each grain $X_i=(X_i,\tau_{X_i})$ \cite{foot3} of which is characterized by its self-similar decomposition space. Namely, all points in $X_i$ are classified into equivalence classes and the equivalence classes coalesce to form a self-similar structure. Then we can obtain a space filled with any $n$ mutually disjoint self-similar spaces. 
Some easily verified and well known statements are shown in the literatures \cite{Nadler}.

Since $(X,\tau)=(\{0,1\}^\Lambda ,\tau_0^\Lambda ),\tau_0=2^{\{0,1\}},Card 
\Lambda \succ \aleph _0$ is easily verified to be a 0-dim, perfect, compact T$_2$-space, there exists a partition \cite{2-a} $\{(X_1,\tau_{X_1}),\dots,(X_n,\tau_{X_n})\}$ of $(X,\tau)$ where each $(X_i,\tau_{X_i})$ is a 0-dim, perfect, compact T$_2$-space \cite{foot4}. 
Then there exists a continuous map $f_i$ from $(X_i,\tau_{X_i})$ onto any compact, self-similar metric space $(Y,\tau_d)$ \cite{2-b}. 
Since $(X_i,\tau_{X_i})$ is a compact space and $Y$ is a T$_2$-space, the map $f_i:(X_i,\tau_{X_i})\rightarrow (Y,\tau_d)$ is a quotient map. 
Therefore the map $h:(Y,\tau_d)\rightarrow (\mathcal{D}_{f_i},\tau(\mathcal{D}_{f_i})),y\mapsto f_i^{-1}(y)$ must be a homeomorphism \cite{2-c}.
Since $(Y,\tau_d)$ is self-similar, the decomposition space $(\mathcal{D}_{f_i},\tau(\mathcal{D}_{f_i}))$ of $(X_i,\tau_{X_i})$ is also self-similar \cite{2-d}.
The space $X$ is filled with self-similar spaces in the sense that the family $\{f^{-1}_i(y);y\in Y,i=1,\dots,n\}$ of subsets of $X$ is a cover of $X$. The self-similar spaces $\mathcal{D}_{f_i}$ and $\mathcal{D}_{f_j}, i\not =j$ are disjoint in the sense that any point $f^{-1}_i(y)$ of $\mathcal{D}_{f_i}$ is not a point of $\mathcal{D}_{f_j}$  as shown in Fig. 1. Here, we note that the decomposition space $\mathcal{D}_{f_i}$ of $X_i$ is not a trivial one $\{\{x\};x\in X_i\}$ especially for a self-similar, connected space $Y$. In fact, disconnected space $\{\{x\};x\in X_i\}$ (a decomposition space of $X_i$) is never homeomorphic to a connected space.

Regarding each subspace $(X_i,\tau_{X_i}),i=1,\dots,n$ as a grain we obtain an aggregate of grains $\{X_1,\dots,X_n\}$ composed of $X_i$ which has a self-similar decomposition space.
Finally, we note that in the above discussions we can replace the self-similar space with a compact substance in the materials science such as a dendrite (A metric space is called a dendrite provided that it is a connected, locally connected, compact metric space which contains no simple closed curve as its subspace. A space which is homeomorphic to a dendrite is a dendrite)\cite{ kitada-csf, Nadler,2-d}
 and then, we can obtain an aggregate of grains $\{X_1,\dots,X_n\}$ composed of grain $X_i$ whose decomposition space is characterized as a dendrite. Fig.1 is available also for the dendritic structure.

The mathematical procedure which plays a central role in the above discussions is summarized as follows. i) $(X,\tau)$ is a 0-dim, perfect, compact T$_2$-space which need not to be metrizable. For example $(X,\tau)=(\{0,1\}^\Lambda ,\tau_0^\Lambda )$ where $Card\Lambda \succ \aleph $ and $\tau_0$ is a discrete topology for $\{0,1\}$. ii) There can exist any $n$ mutually disjoint 0-dim, perfect, compact T$_2$-subspaces $\{(X_1,\tau_{X_1}),\dots,(X_n,\tau_{X_n})\}$ of $(X,\tau)$ and there exist continuous maps $f_i$ from $(X_i,\tau_{X_i}),i=1,\dots,n$ onto any compact metric space $(Y,\tau_d)$ which can be the self-similar structure, e.g., dendrite, each of whose structures are invariant under any homeomorphism. iii) Each decomposition space $(\mathcal{D}_{f_i},\tau(\mathcal{D}_{f_i}))=(\{f^{-1}_i(y)\subset X_i;y\in Y\},\tau(\{f^{-1}_i(y)\subset X_i;y\in Y\}))$ of $(X_i,\tau_{X_i}),i=1,\dots,n$ is homeomorphic to $(Y,\tau_d)$. Then the subspace $(X_i,\tau_{X_i})$ turns to the decomposition space $(\mathcal{D}_{f_i},\tau(\mathcal{D}_{f_i}))$ which shows the self-similar or dendritic structures.

In the present report, we propose the sufficient condition for the problem that a space can be filled with an arbitrary finite number of grains, each of which is characterized as self-similar, using a concept of point set topology. Let a topological space $(X,\tau)$ be a $\Lambda$-product space $(\{0,1\}^\Lambda ,\tau_0^\Lambda )$ of $(\{0,1\},\tau_0)$. Here $Card \Lambda \succ \aleph _0$ (aleph zero) (i.e, $\Lambda $ is an infinite set) and $\tau_0=\{\{0,1\},\{0\},\{1\},\phi \}$. $(X,\tau)$ has an any $n$-partition $\{X_1,\dots,X_n\}$ where each $X_i$ is a non-empty clopen set of $(X,\tau)$ and has a self-similar decomposition space. If we regard each $X_i$ as a grain whose element is a map $x:\Lambda \rightarrow \{0,1\} (\Lambda ={\bf N}, {\bf R},\dots)$, we can mathematically regard the partition $\{X_1,\dots,X_n\}$ of $X$ as an aggregate of grains. Namely, $(X,\tau)$ is filled with any $n$ mutually disjoint self-similar spaces or dendrites each of which is a decomposition space of a grain. To clear the qualitative aspects of the self-similar structure it is emphasized that the discussions do not depend on the metrizability of the initial space $(X,\tau)$.

\bigskip

\noindent
{\bf Acknowledgment}\\
The authors are grateful to Dr. Y. Yamashita, Dr. H. Ryo, Prof. Emeritus H.Fukaishi at Kagawa university, and Prof. H. Nagahama at Tohoku university for useful suggestions and encouragements.


\begin{thebibliography}{9}
\bibitem{2-a}Any 0-dim, perfect, T$_0$(necessarily T$_2$)-space $(X,\tau)$ has its partition $\{X_1,\dots,X_n\}, X_i\in(\tau\cap \Im)-\{\phi\}$ for any $n$. Here $X_i\cap X_{j}=\phi$ for $i\not =j$ and $X_1\cup \cdot \cdot \cdot \cup X_n=X$. Each subspace $(X_i,\tau_{X_i})$ is a 0-dim, perfect, T$_2$-space.
\bibitem{2-b}Let $X$ be a 0-dim, perfect, compact T$_2$-space. Then, for any compact metric space $Y$, there exists a continuous map $f$ from $X$ onto $Y$.
\bibitem{2-c} If $f:(X,\tau)\rightarrow (Y,\tau')$ is a quotient map, then the map $h:(Y,\tau')\rightarrow (\mathcal{D}_f,\tau(\mathcal{D}_f)), y\mapsto f^{-1}(y)$ is a homeomorphism. Here, the decomposition $\mathcal{D}_f$ of $X$ and the decomposition topology $\tau(\mathcal{D}_f)$ are given by $\mathcal{D}_f=\{f^{-1}(y)\subset X; y\in Y\}$ and $\tau(\mathcal{D}_f)=\{\mathcal{U}\subset \mathcal{D}_f;\bigcup \mathcal{U}\in\tau \}$, respectively. See also ref. 6, p. 144.
\bibitem{2-d}A space which is homeomorphic to a self-similar space is self-similar.
\end{thebibliography}
\end{document}